\documentclass[a4paper, 12pt]{article}

\usepackage{amsmath,amsfonts,bm,amssymb,amsthm}
  \usepackage{paralist}
  \usepackage{graphics} %% add this and next lines if pictures should be in esp format
  \usepackage{epsfig} %For pictures: screened artwork should be set up with an 85 or 100 line screen
 \usepackage[colorlinks=true]{hyperref}
\hypersetup{urlcolor=blue, citecolor=red}
\newtheorem{theorem}{Theorem}[section]

\newtheorem{lemma}[theorem]{Lemma}

\theoremstyle{definition}
\newtheorem{definition}[theorem]{Definition}

\newcommand{\ep}{\varepsilon}

\newcommand{\bbR}{\mathbb R}
\newcommand{\bbT}{\mathbb T}
\newcommand{\bbZ}{\mathbb Z}
\newcommand{\const}{\mathrm{const}}
\newcommand{\Ga}{\Gamma}
\newcommand{\gaeps}{\gamma_\ep}
\newcommand{\Diag}{\mathcal D}
\newcommand{\padi}[2]{\frac{\partial #1}{\partial #2}}
\begin{document}

%% Place the running title of the paper with 40 letters or less in []
 %% and the full title of the paper in { }.
\title{Duck farming on the two-torus: multiple canard cycles in generic
		slow-fast systems}

% Place all authors' names in [ ] shown as running head;
% No more than 40 letters. Leave { } empty
% Please use `and' to connect the last two names if applicable
\author{Ilya V. Schurov\thanks{National Research University~--- Higher School of
Economics. This work is supported in part by RFBR (project
10-01-00739), joint RFBR/CNRS
project 10-01-93115-CNRS-a and Grant of President of Russia MK-2790.2011.1.}}
\date{}
% It is required to enter MSC and Keywords.
% Email address of each of all authors is required.
% You may list email addresses of all other authors, separately.

% Put your short thanks below. For long thanks/acknowlegements,
%please go to the last acknowlegments section.

\maketitle

% The name of the associate editor will be entered by an editorial staff
% "Communicated by the associate editor name" is not needed for special issue.
% \centerline{(Communicated by the associate editor name)}

%The abstract of your paper
\begin{abstract}
Generic slow-fast systems with only one (time-scaling) parameter on the
two-torus have attracting canard cycles for arbitrary small values of this
parameter. This is in drastic contrast with the planar case, where canards
usually occur in two-parametric families. In present work, general case of
nonconvex slow curve with several fold points is considered. The number of
canard cycles in such systems can be effectively computed and is no more than
the number of fold points. This estimate is sharp for every system from some
explicitly constructed open set.
\end{abstract}

\vskip 0.7pc
\hskip 1.3pc {\small 2010 Mathematics Subject Classification: Primary 34E15,
37G15.}

{\small \hskip 1.3pc 
Keywords: slow-fast systems, canards, limit cycles, Poincar\'e map, equation of
 variations
 }

%The title of your section 1
\section{Introduction}
Consider a generic slow-fast system:
\begin{equation}\label{eq:main}
	\begin{cases}
		\dot x=f(x,y,\ep)\\
		\dot y=\ep g(x,y,\ep)\\
	\end{cases}\quad \ep\in(\bbR,0).
\end{equation}
For the planar case (i.e. $(x,y)\in \bbR^2$), there is a rather simple description of its
behavior for small $\ep$. It consists of interchanging phases of slow motion
along stable parts of the slow curve $M:=\{(x,y)\mid f(x,y,0)=0\}$ and fast
jumps along straight lines $y=\const$. (See e.g.~\cite{MR}.) Given additional
parameters, depending
on $\ep$, one can observe more complicated behavior: appearance of \emph{duck}
(or \emph{canard}) solutions (particularly limit cycles), i.e. solutions, whose
phase curves contain an arc of length bounded away from 0 uniformly in $\ep$,
that keeps close to the unstable part of the slow curve~\cite{BCDD}.
.

In~\cite{GI}, Yu.~S.~Ilyashenko and J.~Guckenheimer discovered a new kind of
behavior of slow-fast systems on the two-torus. It was shown that for some
particluar family with no auxiliary parameters there exists a sequence of
intervals accumulating at $0$, such that for any $\ep$ from these intervals, the
system has exactly two limit cycles, both of which are canards, where one is
stable and the other unstable.  Yu.~S.~Ilyashenko and J.~Guckenheimer
conjectured that there exists an open domain in the space of slow-fast systems
on the two-torus with similar properties. Here we prove this conjecture, and
provide almost complete description for bifurcations of canard cycles on the
two-torus. In particular, we give sharp estimate for the number of canard cycles
in such systems.

Our main results are the following ones. Consider the system~\eqref{eq:main} and assume that
the phase space is the two-torus:
\begin{equation}
	(x,y)\in\bbT^2\cong \bbR^2/(2\pi\bbZ^2).
\end{equation}
Assume that the speed of the slow motion is bounded away from zero ($g>0$), the
slow curve $M$ is a smooth connected curve, and its lift to the covering
coordinate plane is contained in the interior of the fundamental square
$\{|x|<\pi,\ |y|<\pi\}$. We also assume that all fold points of the slow curve (i.e. the points of
$M$ where the tangent line to $M$ is parallel to $x$-axis) are nondegenerate
(i.e. the tangency rate is quadratic). In this case, the number of fold points
is finite and even: let us denote it by $2N$.

\begin{theorem}
	\label{thm:main}
	For any generic slow-fast system on the two-torus with the 	properties
	described above, under some additional nondegenericity assumptons, the following
	properties hold. There exists a positive number $k\le N$ and a sequence of
	intervals accumulating to zero, such that for every $\ep$ belonging to
	one of these intervals the system has exactly $k$ attracting and $k$ repelling limit
	cycles. All these cycles are canards, and make exactly one turn along the
	$y$-axis during the period.  The measure of their basins of attraction or
	repulsion is bounded away from 0
	uniformly in $\ep$.  Finally, for any sufficiently small $\ep>0$, the number of limit
	cycles that make one turn along the $y$-axis does not exceed $2k$.
\end{theorem} 
\begin{theorem} \label{thm:sharp} There exists an open set in the space of slow-fast systems on
	the two-torus for which the maximal number $k$ of pairs of canard cycles
	reaches its maximal possible value~$N$.  
\end{theorem}
The paper is organized as follows. In section~\ref{sec:descr-of-phen} we provide
an heuristic description of the phenomena discussed. In section~\ref{sec:poincare}
some preliminary results about the Poincar\'e map are stated. Section~\ref{sec:shape}
gives an overview of the proof of Theorem~\ref{thm:main}. This proof relies on
some auxiliary results, that are discussed in sections~\ref{sec:neutral}
and~\ref{sec:order}. Section~\ref{sec:duckfarm} is devoted to construction of
system with maximal number of canards and proof of Theorem~\ref{thm:sharp}.

Due to size limitations, we omit technical details from presented proofs. We
refer the reader to the works~\cite{Sch,Sch2} for more detailed discussion.

\section{Description of the phenomena}\label{sec:descr-of-phen}
In this section we provide heuristic description of the phenomena discovered by
Ilyashenko and Guckenheimer. In what follows, we will assume that $x$-axis of
fast motion is vertical, and $y$-axis is horizontal. The slow motion is directed
from the left to the right.

\begin{figure}[htp]
	\begin{center}
		\includegraphics[scale=1.4]{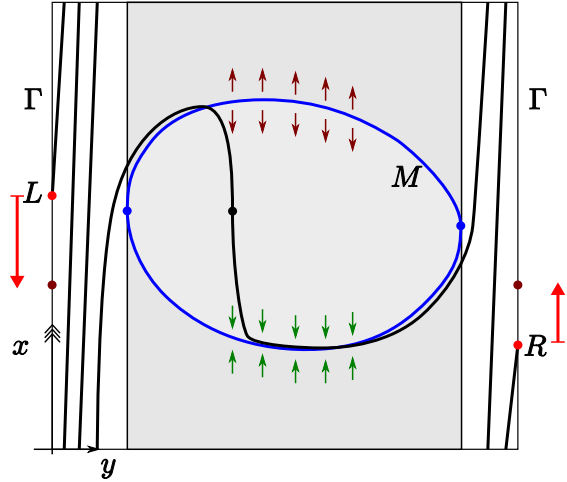}
		\caption{Canard solution of the system with convex slow curve: the base strip
		is shadowed}\label{fig:convex}
	\end{center}
\end{figure}

We consider first the simplest case: $M$ is a \emph{convex curve} and therefore
it has exactly two fold points (i.e. $N=1$).\footnote{In fact, only the
latter condition matters: any system with $N=1$ can be considered as a system
with convex slow curve.} The right one is called
\emph{jump point} and the left one is \emph{reverse jump point}. Consider a
strip $B$ in the phase space that contains $M$ and bounded by vertical circles
that pass through the fold point (see Fig.~\ref{fig:convex}). We call it
\emph{the base strip}. In more generic (nonconvex) case the base strip is defined as
the minimal vertical strip that contains $M$. 

Fix some vertical cross-section $\Ga=\{y=\const\}$ that does not interset $M$.
We will assume without loss of generality that $\Ga=\{y=-\pi\}=\{y=\pi\}$.
Consider some point $w\not\in M$ from the interior of the base strip
$\mathring{B}$.
Trajectory, passing through this point, in forward time is attracted quickly to
the stable part of the slow curve, then moves slowly to the right until 
reaches the jump point, then ``jumps'' and continues slow motion along the $y$-axis,
making about $1/\ep$ rotations along the $x$-axis before it intersects $\Ga$
(call this phase \emph{after-jump rotations}). For given $\ep$, denote the point
of the first intersection with $\Ga$ by $R(\ep)$.

In backward time, the trajectory is quickly attracted to the~\emph{repelling} part of the slow
curve, then moves slowly to the left until reaches the reverse jump point, jumps, and
continues slow motion along the $y$-axis while rotating along the $x$-axis, up to
the intersection with $\Ga$. Denote the point of the intersection by $L(\ep)$. This
trajectory is \emph{canard}, because it has a segment which is close to the repelling
part of the slow curve.

As $\ep>0$ decreases, ``fast'' parts of the trajectory become more vertical, and
the number of rotations during the after-jump motion increases.  Therefore, the point
$R(\ep)$ moves upwards, and $L(\ep)$ moves downwards. By continuity, there
exists $\ep_1$ such that for $\ep=\ep_1$ these two points coincide:
$R(\ep_1)=L(\ep_1)$. This gives us canard limit cycle. As $\ep>0$ continues
decreasing, new coincidense occurs for some $\ep=\ep_2$, $0<\ep_2<\ep_1$, and so
on. Therefore, for the sequence of parameter values $\ep=\ep_k$, accumulating to 
$0$, the system has canard limit cycles. By choosing initial point close enough
to the reverse jump point, it is possible to make these cycles stable. When we
perturb initial point slightly, corresponding values of $\ep_k$ also
perturb slighly, giving us ``canard intervals'', whose existence is stated in
Theorem~\ref{thm:main}.

\begin{figure}[htp]
	\begin{center}
		\includegraphics[scale=0.9]{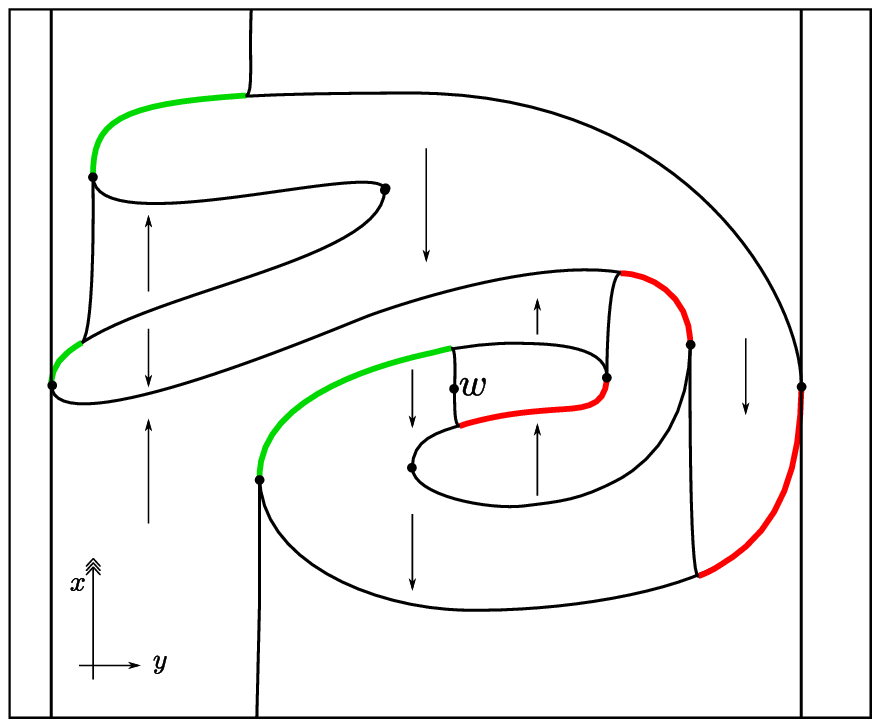}
		\caption{Nonconvex case: several jumps}\label{fig:Z-of-w}
	\end{center}
\end{figure}

When we consider more general case of nonconvex $M$ and $N>1$, the description
becomes more complicated, but the main arguments still work. Let us assume that
any two fold points lie on different vertical circles, and the initial point $w$ does not lie above or below any
fold point. In this case, the trajectory which starts in $w$, in forward
(backward) time falls to attracting (repelling) segment of $M$, moves slowly to
the right (left) until reaches the fold point, jumps and either leaves the base
strip or falls to other attracting (repelling) segment of $M$, and the
process repeats until the trajectory leaves the base strip (see
Fig.~\ref{fig:Z-of-w}).

The main difference with the convex case here is the possibility of several
jumps. In fact, it
does not affect heuristic arguments presented above, because they deal
mostly with the after-jump rotations. However, to provide a rigorous proof of
the main
results and in particular to calculate the number of limit cycles, it does not
suffice to use only the ideas discussed. Instead, we have to perform accurate
analysis of the Poincar\'e map from $\Ga$ to itself, which is discussed in the
next sections.

%In order to calculate the number of limit cycles, we have to analyse the
%Poincar\'e map from $\Ga$ to itself. We discuss it in the next section.

\section{Poincar\'e map}\label{sec:poincare}
Note that the function $g$ is bounded away from zero, so we can divide the
system~\eqref{eq:main} by $g$, thus re-scaling the time: this does not change
the desired properties of its solutions (we are interested only in phase
curves), and the system with the new function $f$ will satisfy the same
nondegenericity assumtions. Thus without loss of generality we can assume $g=1$
in~\eqref{eq:main}.

Consider Poincar\'e map $P_\ep:\Ga\to\Ga$. The slow motion is constant (and bounded
away from $0$), so $P_\ep$ is a well-defined diffeomorphism of a circle.
Its periodic (in particular, fixed) points correspond to closed solutions of the
system.  Denote the graph of $P_\ep$ by $\gaeps$. Fixed points of the Poincar\'e
map correspond to the intersection points of the graph with the diagonal
$\Diag:=\{y=x\}$.  Note, that in terms of the previous section,
$P_\ep(L(\ep))=R(\ep)$.

The derivative of the Poincar\'e map in any point $x_0\in \Ga$ can be easily calculated by
integrating the equation of variations. Namely, let $x=x(y;x_0,\ep)$ be the
phase curve with the initial
condition $x(-\pi;x_0,\ep)=x_0$. Then, 
$P_\ep(x_0)=x(\pi;x_0,\ep)$
and
\begin{equation}\label{eq:P'}
	P'_\ep(x_0)=\left.\padi{P_\ep(x)}{x}\right|_{x=x_0}=\padi{x}{x_0}(\pi;x_0)=\exp\frac{1}{\ep}\int_{-\pi}^{\pi}
	f'_x(x(y;x_0,\ep),y,\ep)\,dy.
\end{equation}
Near the attracting (repelling) parts of the slow curve $M$, the function under the
intergral sign is negative (positive), and the trajectories attract (repell) each other
while moving in these areas. Corresponding parts of the trajectories contribute
contraction (expansion) to the derivative of Poincar\'e map. In ``most of the cases'' 
$\int f'_x(x,y,\ep) dy\ne 0$ and either contraction or expansion dominates,
thus giving either exponentially small or exponentially big derivative with
respect to $\ep$.

It turns out that it is possible to replace actual trajectory in the right-hand
side of~\eqref{eq:P'} with so-called \emph{singular trajectory} (or
\emph{contour}), which is defined as follows. For every point $w\in B\setminus
M$, which does not lie above or below any fold point of $M$, recall the
description of the trajectory which pass through $w$ traced in backward and
forward time up to exit from $B$ (see section~\ref{sec:descr-of-phen}). Assume
that all phases of fast motion in this description are strictly vertical. Then
we obtain a picewise-smooth curve in the base strip which consists of vertical
segments and arcs of the slow curve M, interchanging each other. Call this curve
singular trajectory (or contour) of $w$ and denote it by $Z(w)$. This curve is
in a sense a limit (as $\ep\to 0$) of the trajectories with the initial condition $w$.
The part of the contour to the right of $w$ (which corresponds to the trajectory in
forward time) is denoted by $Z^-(w)$, and the part to the left of $w$ (which
corresponds to backward time) is denoted by $Z^+(w)$.

The following Lemma represents the fact that the derivative of the Poincar\'e map
is controlled (with given precision) by the contour of the corresponding trajectory.
It means that the main contribution to the derivative is made by the segments of slow
motion near the arcs of the slow curve. This contribution dominates over the one
of the jumps
and the after-jump rotations.
\begin{lemma}\label{lem:Q-est}
	Fix some $\delta>0$. Fix some vertical interval $J$, which intersects
	attracting part of $M$ and $\delta$-bounded from repelling part of $M$ (and
	therefore from fold points). Let $u$ be coordinate on $J$. Consider Poincar\'e
	map $Q:J\to \Ga$ in forward time. Then for any $w\in J$,
	\begin{equation}\label{eq:dQdu}
		\log\left.\frac{dQ(u)}{du}\right|_{J}=\frac{1}{\ep}\left[\int_{Z^-(w)}
		f'_x(x,y,0)\, dy+o(1)\right].
	\end{equation}
	Obviously, $Z^-(w)$ does not depend on choice of $w\in J$. The remainder term
	$o(1)\to 0$ as $\ep\to 0$ uniformly in choice of $J$ provided $\delta$
	is fixed. 
\end{lemma}
\begin{lemma}\label{lem:P-est} Let $w\in B\setminus M_\delta$, where $M_\delta$ is
	$\delta$-neighborhood of $M$, and $y(w)$ ($y$-coordinate of $w$) is
	$\delta$-far from $y(G)$ for any fold point $G$. Consider actual trajectory
	that pass through $w$, and denote its initial condition on $\Gamma$ by $x_0$.
	Then
	\begin{equation}
		\log P'_\ep(x_0)=\frac{1}{\ep}\left[\int_{Z(w)}
		f'_x(x,y,0)\,dy+o(1)\right].
	\end{equation}
\end{lemma}
The proof of these Lemmas is given in~\cite{Sch2} (see Lemma 5.4 and Lemma 7.3
there), and rely heavily on the analysis of the after-jump rotations from~\cite{Sch}
(see Theorem 4.3 there, which is reformulated as Theorem 4.6 in~\cite{Sch2}).
%For acknowledgements section, please don't number the section, please begin it with \section*{Acknowledgements}

\section{Shape of the graph of $P_\ep$}\label{sec:shape}
Our goal is to describe the shape of $\gaeps$ and its dependence on $\ep$. The
following description goes back to Shape Lemma in~\cite{GI}, where it is proved
for a particular example.

\begin{figure}[htp]
	\begin{center}
		\includegraphics[scale=0.8]{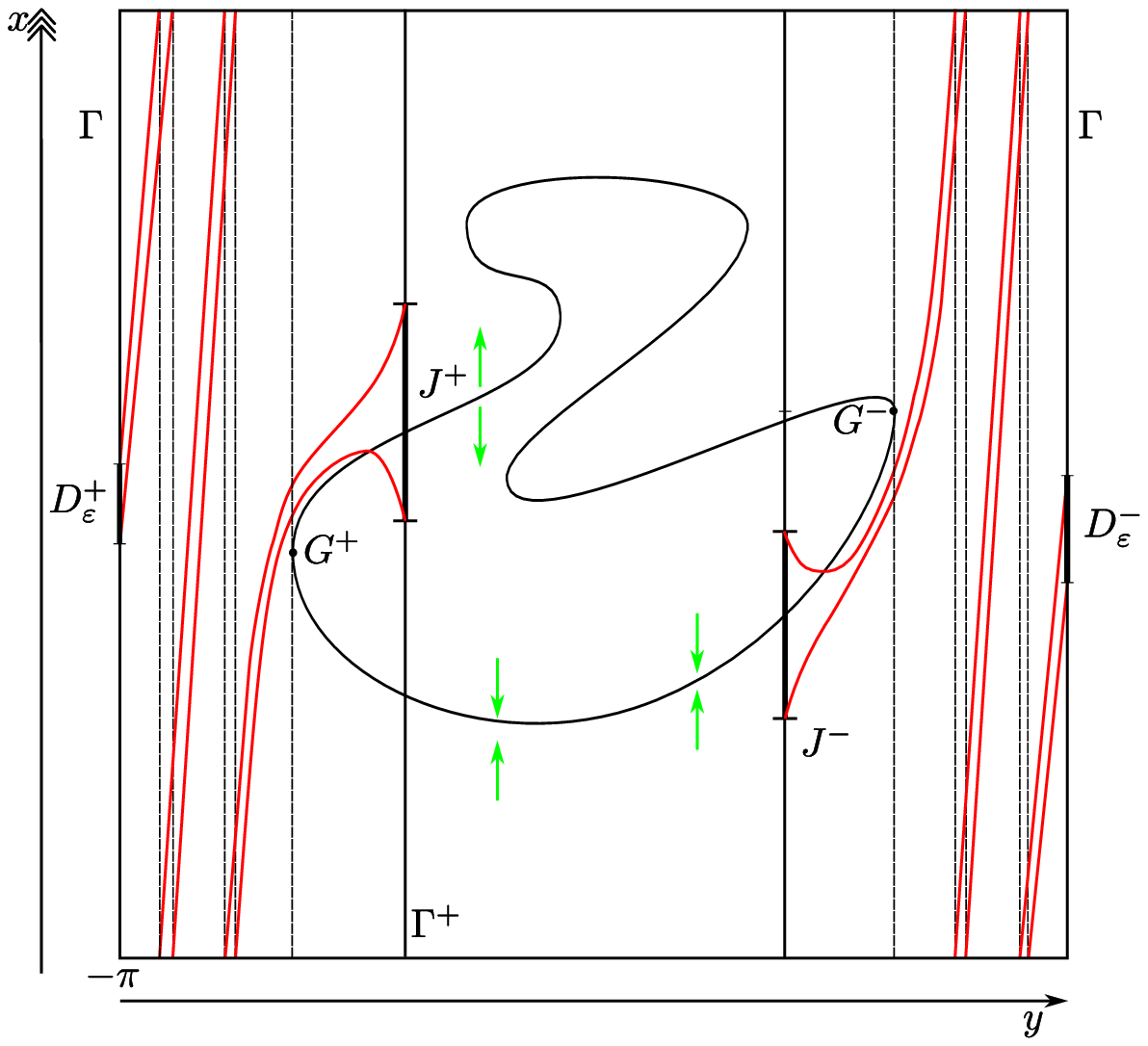}
		\caption{Phase curves near main jump points}\label{fig:general-view}
	\end{center}
\end{figure}

Fix some vertical interval $J^+$ (resp. $J^-$) in the phase space, which
intersects repelling (attracting)
part of the slow curve close enough to far left (far right) fold point and is bounded from
attracting (repelling) part of the slow curve. (See Fig.~\ref{fig:general-view}.)

\begin{definition} A trajectory is called a
	\emph{duck} (\emph{canard}) if and only if it intersects~$J^+$.
\end{definition} 
Consider the projection of $J^+$ (resp., $J^-$) to $\Gamma$ along phase curves
in backward (forward) time.  Denote it by $D^+_\ep$ (resp, $D^-_\ep$). Note that
all the trajectories that intersect $D^+_\ep$ are ducks.  Lemma~\ref{lem:Q-est}
(applied to the system with the time reversed if necessary) implies immediately that $|D^+_\ep|=O(\exp(-C_1/\ep))$ and
$|D^-_\ep|=O(\exp(-C_2/\ep))$ for some positive $C_1,\ C_2$. The trajectory with
the inital condition $x_0\in \Gamma\setminus D^+_\ep$ does not intersect $J^+$.
Therefore, it is attracted to the attracting part of the slow curve rather
quickly (the speed is controlled by the distance between $J^+$ and the reverse jump point), and then
moves near attracting parts of the slow curve only, accumulating contraction
(see the integral~\eqref{eq:P'}). It
also have to intersect $J^-$ and therefore $D^-_\ep$. Lemma~\ref{lem:P-est}
implies that in this case the derivative of the Poincar\'e map is exponentially small.
On the other hand, appliyng the same arguments to the system with time reversed,
we have that outside of $D^-_\ep$, the inverse Poincar\'e map $P^{-1}_\ep$ has
exponentially small derivative. Informally speaking, it means that almost all
the circle in the pre-image (except for very small interval) is mapped into very
small interval in the image, while the exceptional interval in the pre-image is
mapped into almost all the circle in the image.

\begin{figure}[htp]
	\begin{center}
		\includegraphics[scale=1]{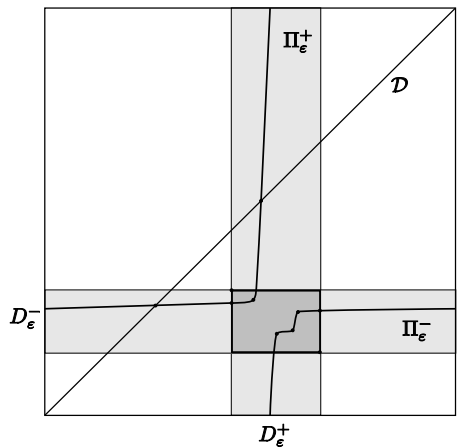}
		\caption{Graph of Poincar\'e map}\label{fig:poincare-graph-zoomout}
	\end{center}
\end{figure}

Geometrically, this means that the graph $\gaeps$ belongs to the union $\Pi^+\cup
\Pi^-$ of exponentially thin strips: vertical
$\Pi^+=D^+_\ep\times S^1$ and horizontal $\Pi^-=S^1\times D^-_\ep$. Outside of
the rectangle $K_\ep=\Pi^+\cap\Pi^-$, the slope of $\gaeps$ is either
exponentially big or exponentially small (see
Fig.~\ref{fig:poincare-graph-zoomout}).

Monotonicity arguments similar to the ones discussed in section~\ref{sec:descr-of-phen}
show that as $\ep\searrow 0$, rectangle $K_\ep$ moves from bottom-right to
top-left corner, making infinitely many rotations.  (See Monotinicity Lemmas
in~\cite{GI} and~\cite{Sch} for details.) 

In this paper, we are interested only in limit cycles that correspond to fixed
points of Poincar\'e map (i.e. making $1$ rotation along the $y$-axis). They are born (or
die) when the diagonal $\Diag$ tangents the graph $\gaeps$. Such a tangency is
possible only at the points where the slope of $\gaeps$ is equal to $1$. We will call
such points \emph{neutral}, applying this term both to the points on the graph $\gaeps$ and
to the corresponding values of argument (i.e. roots of the equation $P'_\ep(x)=1$). Note
that all neutral points belong to $K_\ep$, and therefore the fixed points can be born
only inside $K_\ep$, thus giving us pairs of repelling and attracting canard
cycles. (All points in $K_\ep$ correspond to canard solutions because they lie
over $D^+_\ep$.) 

For every neutral point $x$, consider second derivative $P''_\ep(x)$, and call
$x$ \emph{generating} (resp., \emph{annihilating}) if $P''_\ep(x)<0$ (resp.,
$P''_\ep(x)>0$). We may impose nondegenicity conditions, such that
$P''_\ep(x)\ne 0$ in every neutral point $x$. Consider a projection
$\Delta(x,y)=x-y$ along the diagonal $\Diag$. It turns out (this will be discussed
below) that for any given system for $\ep$ small enough the number of neutral
points is alwas the same (does not depend on $\ep$; see section~\ref{sec:neutral}) and
the order of their projections under $\Delta$ is fixed as well (see
section~\ref{sec:order}). In this case, actual maximal number of canard cycles
is defined by the order of births and deaths, which is controlled by the order
of generating and annihilating neutral points under the projection map $\Delta$, and thus does
not depends on $\ep$. This gives us the number $k$ from Theorem~\ref{thm:main}.
Rolle's Theorem implies that $k\le N$. Neutral points depend on $\ep$
continiously, therefore on every turn of $K_\ep$ there exists an open interval
of $\ep$'s, on which the maximal number (which is $2k$) of canard cycles are born. Such
intervals accumulate to $0$, and their existence is the main result of
theorem~\ref{thm:main}. 

The rest of paper
is devoted to the analysis of neutral points. We first prove that the number
of neutral points is bounded by the number of folds of $M$ and show how it can
be calculated explicitly (see section~\ref{sec:neutral}). Then we discuss the order of births and deaths of
canard cycles (see section~\ref{sec:order}). Finally, we will contstruct an open set of systems with maximal
number of limit cycles $k=N$ (section~\ref{sec:duckfarm}).

\section{Neutral points}\label{sec:neutral}
Consider a trajectory, that passes through some point $w\in\mathring{B}\setminus M$ (see the
description in section~\ref{sec:descr-of-phen}). The part of the trajectory to
the left from $w$ lies near the repelling arcs of the slow curve; the part to the
right from $w$ lies near the attracting parts of the slow curve. We will say
that at
$w$ the trajectory passes through \emph{duck (or canard) jump}: the transition from unstable
part of the slow curve to the stable one. Let
%$S=\{(x,y)\mid x=s(y), y\in [a,b]\}$ 
$S$
be an arc of the slow curve
$M$ between two consequent fold points (\emph{maximal arc}). It is
well-known~\cite{Fenichel} that there exists
invariant curve $S_\ep$ that tends to $S$ as $\ep\to 0$. This curve is called
\emph{(maximal) true slow curve}.
%$S_\ep=\{(x,y)\mid x=s_\ep(y), y\in [a,b]$ which is $O(\ep)$-close to $M'$ in the following sense: for
%every $[a',b']\Subset[a,b]$,  $|s_\ep(y)-s(y)|=O(\ep)$, where $O(\ep)$ is
%estimated uniformly with respect to $y\in [a',b']$. 
It is not unique, but all such curves are exponentially close to each other and
we can pick a suitable one.

\begin{figure}[tp]
	\begin{center}
		\includegraphics[scale=1]{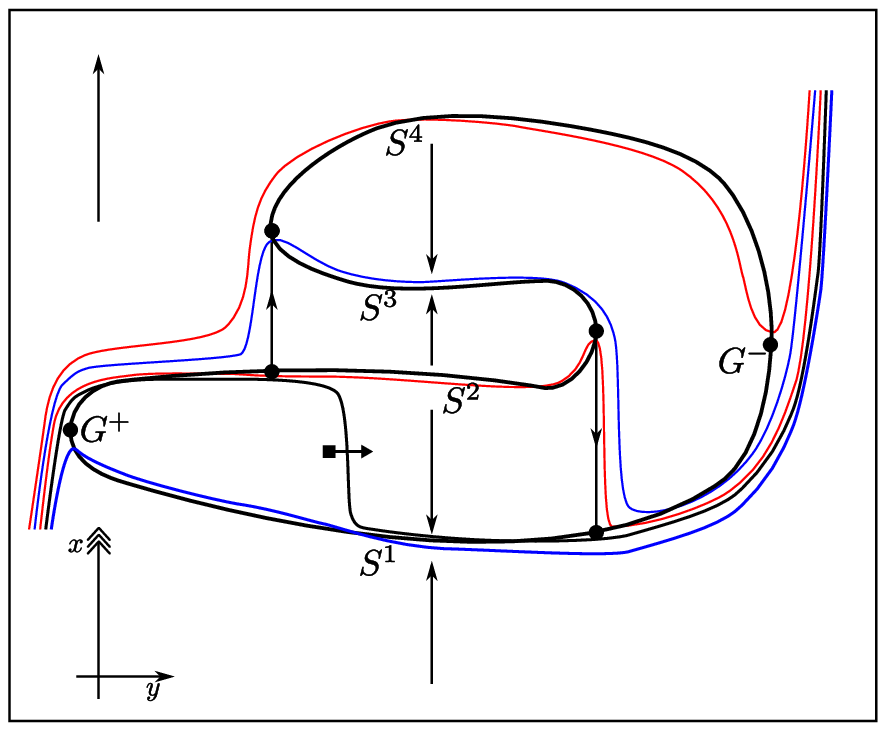}
		\caption{Maximal true slow curves and duck jump}\label{fig:simple-fold-trueslowcurve}
	\end{center}
\end{figure}

For every maximal arc of $M$, consider corresponding maximal true slow curve
(see Fig.~\ref{fig:simple-fold-trueslowcurve}).

Extend them in the backward time to $\Gamma$, and denote corresponding intersection
points by $u_1,\ldots,u_{2N}$ (enumeration is consequent, even numbers correspond
to the repelling curves and odd to the attracting ones; obviously, they should
interchange). Put by definition $u_{2N+i}\equiv u_i$.  Enumerate corresponding slow
curves as $S^1,\ldots, S^{2N}$, and true slow curves as $S^1_\ep,\ldots,
S^{2N}_\ep$ respectively.

The trajectory, that passes through canard jump from $S_{2l}$ can fall after
the jump either to $S_{2l-1}$ or to $S_{2l+1}$. Consider the first case. It becomes 
possible if the initial condition $u$ belongs to the interval $(u_{2l-1},u_{2l})$,
i.e. lies below $u_{2l}$. When we move $u$ a bit upward (closer
to $u_{2l}$), the whole trajectory moves closer to $S^{2l}_\ep$. Thus the
duck jump moves to the right. It means that the trajectory will spend more time
near the repelling part of the slow curve and less time near the attracting part.
Therefore, it will accumulate more expansion and less contraction, and the
derivative of the Poincar\'e map increase monotonically on this
interval.\footnote{To be honest, we are cheating here a little bit: this proves monotonicity only for some smaller
intervals of the vertical circle. Fortunately, they contains all the neutral points, so the proof works.
See~\cite{Sch2} for details.}

Similar arguments show that the derivative of the Poincar\'e map decreases 
monotonically on interval $(u_{2l}, u_{2l+1})$. It follows that the Poincar\'e map
has picewise-monotonic derivative with exactly $N$ intervals of growth and $N$
intervals of decrease. Therefore, the equation $P'(x)=1$ can have no more than $2N$
roots. This proves the estimate for the number of the neutral points, and therefore
of the canard cycles.

Actual number of neutral points can be
calculated as the number of zeroes of logarithmic derivative $\log
P_\ep'(x)$. It follows from the analysis above that this derivative,
presented as an integral~\eqref{eq:P'}, reaches its maximal and minimal values
in the points that corresponds to the maximal true slow curves. Lemma~\ref{lem:P-est} (with some
modifications) implies that these values can be calculated as integrals over special contours, which contain maximal arcs of the
slow curve. Thus the number of the neutral points equals to the number of sign
changes of these values. It does not depend on $\ep$ and can
be effectively calculated.
\section{Order of neutral points}\label{sec:order}
Lift the rectangle $K_\ep$ to the universal cover of the two-torus continuosly
with respect to $\ep$. Pick two arbitrary neutral points $\xi,\,\eta\in \gaeps$
from this lifted rectangle. Then we can define the difference between their
projections $\Delta(\xi)-\Delta(\eta)$, assuming that $\Delta=x-y$, where $x$ and $y$
are coordinates on universal cover. (We need these precuations, because in
general case the difference between two points on a circle is not defined.) In
this section, we show that for any two neutral points the sign of this
difference does not depend on $\ep$.

The main idea is to show that the segment $[\xi, \eta]\subset K_\ep$ is either
``almost vertical'' or ``almost horizontal''. In the first case, if $\xi$ is top end
of the segment and $\eta$ is bottom end, then $\Delta(\xi)-\Delta(\eta)>0$. In the
second case, if $\xi$ is left end and $\eta$ is right end, then
$\Delta(\xi)-\Delta(\eta)>0$, and so on.

Consider two trajectories which correspond to neutral points (call them neutral
trajectories). Due to Lemma~\ref{lem:P-est}, they should lie near some contours with zero
integrals (call such contours neutral as well). Note, that this implies
that the measure of basins of limit cycles is bounded away from $0$: these
cycles lie in different
areas of the phase space, which are separated by neutral solutions which are close to fixed neutral
contours. 

Consider first forward-time parts of
these contours (which are denoted by $Z^+$). They both contain the far right
fold point and therefore have some nonempty intersection. Denote far left
point of this intersection by $T$. To the left of $T$, the corresponding contours (and
therefore actual trajectories) are bounded away from each other. To the right of $T$, the trajectories follow the same attracting arcs of the slow curve, and
therefore attract each other. 

Consider Poincar\'e map from some interval $J'\ni T$
to $\Ga$ in forward time. Then the rate of the attraction is given by
Lemma~\ref{lem:Q-est} and is defined by the integral $f'_x$ over intersection of the
contours. This integral does not depend on $\ep$ and can be calculated
explicitly. The distance between these trajectories when they approach $\Gamma$
is the distance between the $y$-coordinates of corresponding neutral points. It
follows immediately that $|y(\xi)-y(\eta)|=O^*(\exp(-C^-/\ep))$ for some
$C^->0$.

Applying the same arguments to the system with time reversed, we obtain similar
statement for $x$-coordintes: $|x(\xi)-x(\eta)|=O^*(\exp(-C^+/\ep))$ for some
$C^+>0$.

Consider the slope of the segment $[\xi,\eta]$, which is equal to
\begin{equation}
	\frac{|y(\xi)-y(\eta)|}{|x(\xi)-x(\eta)|}=O^*\left(\exp
	\frac{C^+-C^-}{\ep}\right).
\end{equation}
Again, we may impose additional
nondegenericity conditions and assume that $C^+\ne C^-$. This means that the slope is
either exponentially big or exponentially small, what implies the necessary assertion
immediately.

This proves that the order of neutral points under the projection $\Delta$
is fixed and thefore the maximal number of canard cycles $k$ is well-defined. It
finishes the proof of Theorem~\ref{thm:main}.
\section{Duck farm}\label{sec:duckfarm}
The discussion above shows that we can translate ``dynamical'' questions (e.g.
about limit cycles, Poincar\'e map and so on) into
ge\-o\-met\-ri\-cal/com\-bi\-na\-to\-ri\-al
language which involves the shape of the slow curve $M$ and values of integrals
of $f'_x$ over some arcs of $M$. As an application of this approach, we pick
an arbitrary $N>1$ and construct a system with maximal number of canard cycles:
$k=N$. In fact, this example provides an open set of such systems, because all
conditions imposed on the system during the construction of this example are open.  To simplify the
notation, we consider only the case $N=3$, but extension of these arguments to
the general
case is strightforward.

\begin{figure}[tp]
	\begin{center}
		\includegraphics[scale=0.6]{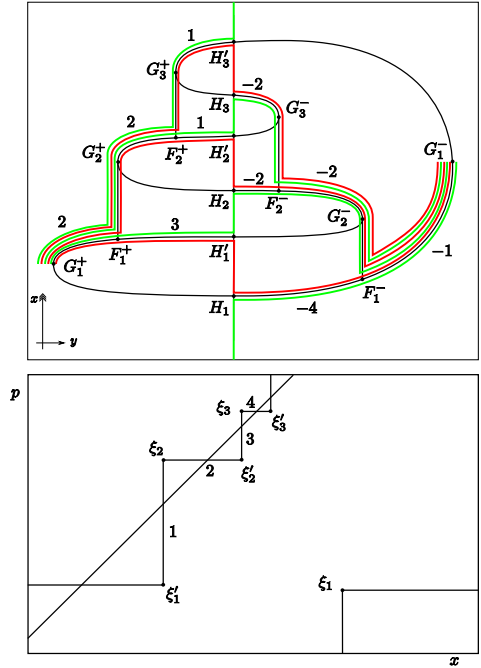}
		\caption{The system with maximal number of
		ducks}\label{fig:maxduck}
	\end{center}
\end{figure}

The key ingredient of the construction is the shape of the slow curve, see
Fig.~\ref{fig:maxduck}, top part. We demand here the depicted contours to be neutral, and
the integrals of $f'_x$ over the corresponding arcs to be equal to corresponding values
(e.g. $\int_{H_1 F_1^-} f'(x,y,0)\,dy=-4$, and so on).

This system has $2N$ neutral contours, and therefore $2N$ neutral points on the
graph $\gaeps$. It follows from the previous results, that for such a system, the
graph looks like a ``staircase'', where ``lengths'' and ``heights'' of the steps
monotonically decrease (see Fig.~\ref{fig:maxduck}, bottom part). This can be
shown by an explicit calculation of the corresponding exponential rates that control
``lengths'' and ``heights'' of the steps (see the description in the previous
section). They depend only on the integrals over the
arcs which we control.

Due to the shape of the graph of Poincar\'e map, it follows that the order of
the bifurcations of limit cycles is the following: first we have $N$ births and then
we have $N$ deaths. During every birth a pair of cycles appear, therefore the
number of canard cycles here is maximal and equal to $2N$. Thus we have constructed
the desired example.  This proves Theorem~\ref{thm:sharp}.

\section*{Acknowledgments} 
The author would like to express his sincere appreciation to Yu.~S.~Ilyashenko
for the statement of the problem and his assistance with the work and to
V.~Kleptsyn for fruitful discussions, valuable ideas and comments both on
mathematics and English of this paper.

% You may incorporate your references as follows in your main tex file.
% Using BibTex is not recommended but can be handled.

\end{document}